\documentclass[12pt,twoside]{article}
\RequirePackage[OT1]{fontenc}
\usepackage{amsthm,amsmath,natbib,amssymb,enumitem,authblk,geometry,amsfonts}
\RequirePackage{hyperref}
 \usepackage{graphicx}
 
 \include {epsf}

\usepackage[dvipsnames]{xcolor}

\newtheorem{theorem}{Theorem}
\newtheorem{corollary}[theorem]{Corollary}
\newtheorem{proposition}[theorem]{Proposition}

\newtheorem{example}[theorem]{Example}

\newcommand{\re}{\mathbb{R}}

\newcommand{\ppp}{\mathcal{P}}
\newcommand{\yyy}{\mathcal{Y}}
\newcommand{\bbb}{\mathcal{B}}

\newcommand{\ponto}{\, \cdot \,}

\newcommand{\ie}{{\it i.e.}}
\newcommand{\deffeq}{\mathrel{\overset{\makebox[0pt]{\mbox{\normalfont\tiny\sffamily def}}}{=}}}

\newcommand{\Feight} {\fontsize{8}{11}\selectfont  }

\begin{document}

\title{Construction and Extension of  Dispersion Models}
 \author{
   Rodrigo Labouriau
   \thanks{Department of Mathematics,
                Aarhus University.\\
                e-mail: rodrigo.labouriau@math.au.dk and rodrigo.labouriau@rlstatlab.com 
           }
  }
\date{September 2026}

\maketitle

\thispagestyle{empty}

\begin{abstract}
\noindent
\Feight
There are two main classes of dispersion models studied in the literature: proper dispersion models (PDMs) and exponential dispersion models (EDMs). Dispersion models that are neither proper nor exponential dispersion models are termed here non-standard dispersion models (NSDMs). This paper exposes techniques for constructing new proper dispersion models and non-standard dispersion models. These constructions provide a solution to an open question in the theory of dispersion models concerning the construction and extension of non-standard dispersion models.

Given a unit deviance function, a dispersion model is usually constructed by calculating a normalising function that makes the density function integrate to one. This calculation involves the solution of non-trivial integral equations. The main idea explored here is to use characteristic functions of real symmetric absolutely continuous probability measures to construct large families of unit deviances for which the associated integral equations admit a trivial solution, in the sense that the normalising function is constant and independent of the observed values. The corresponding dispersion models are proper dispersion models. We then show that non-trivial solutions can also be constructed from spectral properties of the exponential kernel associated with the unit deviance. In particular, zeros of its Fourier transform yield bounded non-constant normalising functions, while zero sets of positive Lebesgue measure give rise to non-trivial $L^2(\mathbb{R})$ orthogonal perturbations. The resulting dispersion models are non-standard models. These constructions show that the classes of proper and non-standard
dispersion models (PDM and NSDM) are vast, with cardinality at least as large as $c=\vert \mathbb{R} \vert$; in particular, the construction of PDMs
contains an infinite-dimensional family parametrised by symmetric
absolutely continuous probability distributions, while explicit
continuum-sized families of NSDMs are also obtained.
\end{abstract}

\noindent
\Feight
\normalsize

\newpage
\tableofcontents
\newpage

\section{Introduction}

Dispersion models (DMs) are parametric families of probability measures defined on Euclidean spaces, which play an essential role in statistics \citep{Jorgensen1987A,Jorgensen1997}.
Several classic families of probability measures such as the normal, Poisson, binomial, gamma, inverse-Gaussian, von Mises, gamma-Poisson compound, and simplex families are DMs. The mathematical properties common to all DMs allow constructing a sound and well-elaborated theoretical machinery to perform statistical inference for those parametric families \citep{Jorgensen1987A,Jorgensen1997,Cordeiroetal2019}.
DMs naturally appear in many statistical applications since they form the basis of several major classes of statistical models such as generalised linear models \citep{McCullagh-Nelder1989,Jorgensen1987A,Jorgensen1987B}, generalised additive models, and variants of state-space models \citep{Jorgensen-etal1996A,Fahmeir2001}.
In this paper, we introduce some techniques for constructing new DMs
which allow us to exhibit large, and in particular infinite-dimensional,
families of DMs; we also obtain lower bounds on the cardinality of some
important classes of DMs.

There are two major classes of DMs studied in detail in the literature: exponential dispersion models (EDMs) and proper dispersion models (PDMs).
Dispersion models that are neither a PDM nor an EDM are termed here \emph{non-standard dispersion models} (NSDMs).
These models ({\it sic}.)
``are still not well understood, mainly for lack of examples of this kind''
\citep[p. 8, last paragraph]{Jorgensen1997}.
Here, we develop general constructions for both PDMs and NSDMs. Characteristic functions of symmetric absolutely continuous distributions are used to generate large classes of unit deviances and PDMs. We also show that non-constant normalising functions, and hence NSDMs, can be constructed from spectral properties of the exponential kernel associated with a unit deviance. In particular, isolated zeros of the Fourier transform yield bounded non-constant normalising functions, whereas zero sets of positive Lebesgue measure give rise to non-trivial $L^2(\mathbb{R})$ orthogonal complements and infinite-dimensional families of perturbations. These constructions imply that the classes of PDMs and NSDMs are large; in particular, continuum-sized families of such models can be obtained explicitly. Some of the ideas developed here were sketched in \citet{Cordeiroetal2019}, but the constructions and the spectral arguments presented below are different and considerably more explicit.

The paper is organised as follows. Section \ref{Sect.2} gives a short review of dispersion models and introduces the basic notation. Section \ref{Sect.3} considers the general problem of constructing a dispersion model from a given unit deviance. Subsection 3.1 discusses the construction of PDMs through the normalising equation. Subsection \ref{Sect.3.2} introduces a general method for generating PDMs from characteristic functions by preserving their local deviance structure while modifying the tails to ensure integrability. Subsection \ref{Sect.3.3} develops a construction of NSDMs based on isolated zeros of the Fourier transform of the exponential kernel, while Subsection \ref{Sect.3.4} studies an $L^2(\mathbb{R})$ version of the construction in which the orthogonal complement of the translation space is characterised by the zero set of the Fourier transform. Section \ref{Sect.4} discusses regularity and illustrates the different constructions by examples.

\section{A Short Review on One-Dimensional Dispersion Models}\label{Sect.2}
Consider a parametric family of  probability measures  $\ppp = \left\{P_{\mu\lambda} : \mu\in\Omega, \lambda\in\Lambda \right\}$ defined on $( \re , \bbb (\re ) )$, where the parametrisation given by  $(\mu,\lambda )$ is identifiable (\ie, the mapping $(\mu,\lambda ) \mapsto P_{\mu\lambda}$  is a bijection between $\Omega\times\Lambda$ and $\ppp$).
Assume that $\ppp$ is dominated by a $\sigma$-finite measure $\upsilon$ with support $\yyy\subseteq\re$, $\Omega \subseteq \yyy$ is open, $\Lambda\subseteq\re_+$ is an interval bounded from the left and unbounded to the right, and that for each $(\mu,\lambda)\in\Omega\times\Lambda$ a version of the Radon-Nykodyn derivative of $P_{\mu\lambda}$ with respect to $\upsilon$ is of the form
\begin{equation}\label{Eq.Oned-01}
\frac{d P_{\mu\lambda}}{d\upsilon} (y)
\deffeq
p(y;\mu,\lambda) = 
a(y;\lambda)  \exp \left \{ - \lambda \, d(y;\mu ) \right \}, 
\mbox{ for all } y\in \yyy.
\end{equation}
Here $a: \yyy\times\Lambda \rightarrow \re_+ $ and $d: \yyy\times\Omega \rightarrow \re_+$ are given suitable functions.
If the function $d$ is such that $d(\mu;\mu) = 0$ for all $\mu\in\Omega$ and $d(y;\mu) >0$ for all $(y,\mu)\in \yyy\times\Omega$ such that $y\ne\mu$, then $d$ is said to be a \emph{unit deviance} and the family $\ppp$ is the 
\emph{dispersion model} (DM) generated by the unit deviance $d$. 
The parameters $\mu$ and $\lambda$ are called the \emph{position parameter} and the \emph{index parameter}, respectively, $\Lambda$ is the \emph{index set}, and the function $a$ is termed the \emph{normalizing function}. Without loss of generality we assume that $\Lambda = \re_+$.

A dispersion model with density (\ref{Eq.Oned-01}) is said to be a \emph{proper dispersion model} generated by a unit deviance $d$ when the normalizing function $a$ factorizes as,
\begin{equation}\nonumber
a(y;\lambda) = a_0(\lambda)\,b(y),
 \mbox{ for all } (y, \lambda) \in \yyy \times \Lambda \, ,
\end{equation}
where $a_0:\Lambda \rightarrow \re_+$ and $b:\yyy  \rightarrow \re_+$ are 
 suitable functions. 
 A dispersion model generated by a  unit deviance  $d$ is said to be an \emph{exponential dispersion model} when the unit deviance takes the form
\begin{equation}\label{eq.2.1.02}
d(y;\mu) = y f(\mu) + g(\mu) + h(y) \, , \mbox{ for all } y\in \yyy \mbox{ and } \mu\in\Omega
\, ,
\end{equation}
for suitable functions $f,g$ and $h$. 
Examples of PDMs are the von Mises, the simplex, the normal, the gamma and the inverse Gaussian families of distributions. The normal, gamma, inverse Gaussian, Poisson and gamma compound Poisson families of distributions are classic examples of EDMs. There are only three PDMs that are also EDMs: the normal, the gamma and the inverse Gaussian families of distributions \citep[Theorem 5.6.]{Jorgensen1997}.

\section{The Problem of Construction of  Dispersion \\ Models}
\label{Sect.3}

We shall be concerned below with the general problem of constructing DMs in a process that will enable us to access the extension of this class of parametric families. Not all unit deviances generate a DM. Indeed, a unit deviance $d$ generates a DM if, and only if, it is possible to find a normalising function  $a:\yyy\times\Lambda \rightarrow \re_+ $  such that the integral of the density given by (\ref{Eq.Oned-01}) integrates $1$, \ie, the function $a$ is the solution of the integral equation
\begin{equation}\label{Eq.Oned-02}
\int_\yyy a(y; \lambda) 
 \exp \left \{ - \lambda d(y;\mu ) \right \} 
 \upsilon (dy) = 1, 
 \mbox{ for all }
 (\mu,\lambda)\in\Omega\times\Lambda \, .
\end{equation}
Note that the solution $a$ should be a function independent of the position parameter $\mu$.
Posed in this generality, this problem is hard to solve since it is difficult to establish even whether the integral equation (\ref{Eq.Oned-02}) has a solution. However, we will present a technique for constructing many examples where the problem of generating dispersion models tractable. 

\subsection{Generation of Proper Dispersion Models}\label{Sect.3.1}

We turn now to the problem of generating proper dispersion models from a given unit deviance $d$. In that case, the normalising function factorises as $a(y;\lambda) = a_0 (\lambda) \, b(y)$, for all $(y, \lambda) \in \yyy\times\Lambda$. Since $p(\ponto; \mu, \lambda)$ is a probability density and  $\exp \left\{ - \lambda d(\ponto ;\mu)  \right\}$ take only positive values, then the function $a = a_0\, . b$ takes only positive values in the support of $\upsilon$. We assume then, without loss of generality, that both $a_0$ and $b$ take only positive values. When working with PDMs the integral equation  (\ref{Eq.Oned-02}) takes the form
\begin{equation}  \label{Eq.Oned-02DM}
  a_0(\lambda) \int_\yyy b(y)   \exp \left\{ - \lambda d(y ;\mu) \right\} \upsilon (dy) =
  1,
 \mbox{ for all } (\mu, \lambda) \in \Omega\times\Lambda
 \, .
\end{equation}
Any function $b:\yyy \longrightarrow \re_+$ such that the integral $\int_\yyy b(y)  \exp \left\{ - \lambda d(y ;\mu) \right\} \upsilon (dy)$ is finite and does not depend on the position parameter $\mu$ generates a proper dispersion model. In the case such a function $b$ exists, defining 
\begin{equation}\nonumber
a_0(\lambda) = \frac{1}{\int_\yyy b(y)  \exp \left\{ - \lambda d(y ;\mu) \right\} \upsilon (dy)}\, , \, \, 
\mbox{ for each $\lambda\in\Lambda$,}
\end{equation}
yields a solution for equation (\ref{Eq.Oned-02DM})  and the densities of the form given by (\ref{Eq.Oned-01}) defined with  $a_0$ and $b$ correspond to the densities of a proper dispersion model.
Note that $\int_\yyy b(y) \exp \left\{ - \lambda d(y ;\mu) \right\} \upsilon (dy) > 0$ since both $ b$ and $\exp \left\{ - \lambda d(\ponto ;\mu) \right\}$ take only positive values in the support of $\upsilon$.
In particular, any function $b\in L^1_\upsilon  (\re)$ such that $\int_\yyy b(y) \exp \left\{ - \lambda d(y ;\mu) \right\} \upsilon (dy)$ does not depend on $\mu$ can be used to construct a PDM.
We present next a technique for obtaining solutions of the equation (\ref{Eq.Oned-02DM}) for a rich class of unit deviances.

\subsection{Generation of Dispersion Models via Unit Deviances Constructed with Characteristic Functions}\label{Sect.3.2}

We present in this section a general construction of proper dispersion
models based on characteristic functions. The main idea is to use a
characteristic function to specify the behaviour of the unit deviance in
a neighbourhood of the diagonal $y=\mu$, and to modify its behaviour
away from the diagonal in such a way that the corresponding exponential
kernel is integrable.

Let $\varphi$ be the characteristic function of an absolutely continuous
probability measure on $\mathbb{R}$ which is symmetric about zero. Then
$\varphi$ is real-valued, continuous and even, $\varphi(0)=1$, and
$ \lim_{|t|\longrightarrow\infty}\varphi(t)=0 .$
Moreover,
\[
        |\varphi(t)|<1, \qquad t\neq0.
\]
Indeed, equality $|\varphi(t_0)|=1$ for some $t_0\neq0$ would imply that
the corresponding probability measure is concentrated on a lattice,
which is impossible for an absolutely continuous probability measure \citep[p. 2, theorem 1.1.4]{Ushakov1999}.
Consequently, $1-\varphi(t)>0$,  for all $t\neq0$, and $1-\varphi(\,\cdot\, )$ is a unit deviance.

For $0<a<b$, let $\rho_{a,b}:\mathbb{R}\longrightarrow[0,1]$ be an
even continuous function satisfying
\[
 \rho_{a,b}(t)=
 \begin{cases}
  0, & |t|\leq a,\\
  1, & |t|\geq b.
 \end{cases}
\]
The function $\rho_{a,b}$ may, for instance, be chosen to be infinitely
differentiable, with a smooth transition between $0$ and $1$ on
$a<|t|<b$.

For each $c>0$, define
\begin{equation}
\label{eq:CFunitdeviance}
 d_{\varphi,a,b,c}(y;\mu)
 =
 1-\varphi(y-\mu)
 +
 c(y-\mu)^2\rho_{a,b}(y-\mu),
 \qquad (y,\mu)\in\mathbb{R}^2 .
\end{equation}
The following proposition shows that this construction generates a
large class of proper dispersion models.

\begin{proposition}
\label{prop:CFPDM}
Let $\varphi$ be the characteristic function of an absolutely continuous
probability measure on $\mathbb{R}$ which is symmetric about zero, and
let $0<a<b$ and $c>0$. Then the function
$d_{\varphi,a,b,c}$ defined by \eqref{eq:CFunitdeviance} is a unit
deviance. Moreover, for every $\lambda>0$,
\[
 \int_{\mathbb{R}}
 \exp\{-\lambda d_{\varphi,a,b,c}(y;\mu)\}\,dy
 <\infty ,
\]
and the value of this integral does not depend on $\mu$. Consequently,
$d_{\varphi,a,b,c}$ generates a proper dispersion model on
$(\mathbb{R},\mathcal{B}(\mathbb{R}))$, dominated by the Lebesgue measure,
with $b(y)=1$.
\end{proposition}

\begin{proof}
Since $\varphi(0)=1$ and $\rho_{a,b}(0)=0$, $d_{\varphi,a,b,c}(\mu;\mu)=0$.
If $y\neq\mu$, then $1-\varphi(y-\mu)>0$,
and the second term on the right-hand side of
\eqref{eq:CFunitdeviance} is non-negative. 
Hence $d_{\varphi,a,b,c}(y;\mu)>0$ for $y\neq\mu$,
and $d_{\varphi,a,b,c}$ is therefore a unit deviance.

Put
$
 d_{\varphi,a,b,c}^0(t)
 =
 1-\varphi(t)+ct^2\rho_{a,b}(t).
$
The unit deviance is of location type,
$
        d_{\varphi,a,b,c}(y;\mu)
        =
        d_{\varphi,a,b,c}^0(y-\mu).
$
It follows, by the change of variable $z=y-\mu$, that
\[
 \int_{\mathbb{R}} \exp\{-\lambda d_{\varphi,a,b,c}(y;\mu)\}\,dy
  =
 \int_{\mathbb{R}} \exp\{-\lambda d_{\varphi,a,b,c}^0(z)\}\,dz ,
\]
which does not depend on $\mu$.

It remains to show that this integral is finite. Since
$\rho_{a,b}(z)=1$ for $|z|\geq b$ and
$1-\varphi(z)\geq0$, we have
$
        d_{\varphi,a,b,c}^0(z)\geq cz^2,
        \qquad |z|\geq b.
$
Therefore, $\exp\{-\lambda d_{\varphi,a,b,c}^0(z)\} \leq \exp\{-\lambda c z^2\}$,  
for $ |z|\geq b$.
The integral over $\{|z|\geq b\}$ is consequently finite, while the
integral over the compact interval $[-b,b]$ is finite because the
integrand is continuous and bounded. 
Thus $ 0< \int_{\mathbb{R}} \exp\{-\lambda d_{\varphi,a,b,c}^0(z)\}\,dz <\infty$.
Defining
\begin{equation}
\label{eq:a0CF}
 a_0(\lambda)
 =
 \left[
 \int_{\mathbb{R}}
 \exp\{-\lambda d_{\varphi,a,b,c}^0(z)\}\,dz
 \right]^{-1},
\end{equation}
we obtain
\[
 p(y;\mu,\lambda)
 =
 a_0(\lambda)
 \exp\{-\lambda d_{\varphi,a,b,c}(y;\mu)\},
\]
which integrates to one for every $(\mu,\lambda)\in
\mathbb{R}\times\mathbb{R}_+$. Hence the corresponding dispersion model
is proper, with $b(y)=1$.
\end{proof}

An important feature of the construction above is that the unit
deviance coincides locally with the function generated directly by the
characteristic function. Namely,
\[
        d_{\varphi,a,b,c}(y;\mu)
        =
        1-\varphi(y-\mu)
\]
whenever $|y-\mu|\leq a$. Thus the additional term in
\eqref{eq:CFunitdeviance} does not alter the local behaviour of the
deviance around the diagonal; its role is only to impose sufficient
growth at infinity to make the exponential kernel integrable.
In particular, if the probability measure associated with $\varphi$
has a finite and non-zero second moment, then
\[
        \varphi(t)
        =
        1-\frac{\sigma^2}{2}t^2+o(t^2),
        \qquad t\longrightarrow0,
\]
where $\sigma^2$ denotes its variance. Hence
\[
        d_{\varphi,a,b,c}^0(t)
        =
        \frac{\sigma^2}{2}t^2+o(t^2),
        \qquad t\longrightarrow0.
\]
Thus the construction preserves the usual quadratic local behaviour of
a regular unit deviance.

The construction also gives a direct lower bound for the cardinality of
the class of proper dispersion models. Let $\mathcal{A}_s$ denote the
class of all absolutely continuous probability measures on
$\mathbb{R}$ that are symmetric about zero. Fix $0<a<b$ and a function
$\rho_{a,b}$ as above. For every $P\in\mathcal{A}_s$, let
$\varphi_P$ denote its characteristic function. Proposition
\ref{prop:CFPDM} associates, for every $c>0$, a PDM to the pair $(P,c)$.

\begin{corollary}
\label{cor:cardinalityPDM}
For each $P\in\mathcal{A}_s$, the construction
\eqref{eq:CFunitdeviance} generates a family of distinct proper
dispersion models having cardinality at least that of $\mathbb{R}$.
Consequently,
\[
        |\mathop{\rm PDM}|
        \geq
        |\mathcal{A}_s|\,|\mathbb{R}|.
\]
In particular, since
\[
        |\mathcal{A}_s|=|\mathbb{R}|,
\]
the class of proper dispersion models has cardinality at least that of
the real line.
\end{corollary}

\begin{proof}
We first show that different pairs $(P,c)$ give different dispersion
models. Suppose that $P_1,P_2\in\mathcal{A}_s$, with characteristic
functions $\varphi_1$ and $\varphi_2$, and let $c_1,c_2>0$.
Write $d_i^0(t) = 1-\varphi_i(t)+c_i t^2\rho_{a,b}(t)$, for $i=1,2$.
Suppose that, for some $\lambda_1,\lambda_2>0$, the corresponding
densities with position parameter zero coincide:
\[
 a_{01}(\lambda_1)
 \exp\{-\lambda_1d_1^0(t)\}
 =
 a_{02}(\lambda_2)
 \exp\{-\lambda_2d_2^0(t)\},
 \qquad \forall t\in\mathbb{R}.
\]
Evaluating this equality at $t=0$ gives $a_{01}(\lambda_1)=a_{02}(\lambda_2)$,
and consequently
\[
        \lambda_1d_1^0(t)=\lambda_2d_2^0(t), 
        \qquad \forall t\in\mathbb{R}.
\]
For $|t|\geq b$,
\[
 \lambda_1\{1-\varphi_1(t)+c_1t^2\}
 =
 \lambda_2\{1-\varphi_2(t)+c_2t^2\}.
\]
Dividing by $t^2$ and letting $|t|\longrightarrow\infty$ yields
\[
        \lambda_1c_1=\lambda_2c_2.
\]
Using this equality in the preceding identity and recalling that
$\varphi_i(t)\longrightarrow0$, we obtain, by letting
$|t|\longrightarrow\infty$ once more,
\[
        \lambda_1=\lambda_2.
\]
Hence $c_1=c_2$, and the equality
$\lambda_1d_1^0=\lambda_2d_2^0$, and we obtain that
\[
        \varphi_1(t)=\varphi_2(t), 
        \qquad \forall t\in\mathbb{R}.
\]
By the uniqueness theorem for characteristic functions,
$P_1=P_2$.

Thus the mapping
\[
        (P,c)\longmapsto \hbox{the PDM generated by }
        d_{\varphi_P,a,b,c}
\]
is injective. Since $|\mathbb{R}_+|=|\mathbb{R}|$, every
$P\in\mathcal{A}_s$ generates at least $|\mathbb{R}|$ distinct PDMs,
which proves the first assertion.

Finally, the class $\mathcal{A}_s$ has the cardinality of the real
line. Therefore,
\[
        |\mathop{\rm PDM}|
        \geq
        |\mathcal{A}_s|\,|\mathbb{R}|
        =
        |\mathbb{R}|\,|\mathbb{R}|
        =
        |\mathbb{R}|.
\]
\end{proof}

\noindent
Here are two remarks on the construction above:
First, the choices of the transition interval $(a,b)$ and of the cut-off
function $\rho_{a,b}$ provide additional freedom in the construction.
The parameter $c$ was used in Corollary \ref{cor:cardinalityPDM}
because it gives a particularly simple injective parametrisation of a
continuum of PDMs for each characteristic function.
Second, the quadratic term in \eqref{eq:CFunitdeviance} is not essential.
It may be replaced by any non-negative even function $q$ satisfying
$q(t)=0$ in a neighbourhood of zero and growing sufficiently rapidly
that 
$\exp\{-\lambda q(t)\}$
is integrable over the real line for every $\lambda>0$. The use of a quadratic
term is convenient because it gives an immediate comparison with a
Gaussian kernel.

There is another sense in which the class of PDMs obtained by this
construction is large. Let $\mathcal{A}_s$ denote the class of
absolutely continuous probability measures on $\mathbb{R}$ that are
symmetric about zero. Identifying these measures with their densities,
$\mathcal{A}_s$ is an infinite-dimensional convex subset of
$L^1(\mathbb{R})$; more precisely, its affine hull is
infinite-dimensional. Fixing $a$, $b$, $c$ and $\rho_{a,b}$, the map
\[
 P\longmapsto
 d_P^0(t)
 =
 1-\varphi_P(t)+ct^2\rho_{a,b}(t),
 \qquad P\in\mathcal{A}_s,
\]
where $\varphi_P$ denotes the characteristic function of $P$, is
injective and affine. Indeed, characteristic functions determine
probability measures uniquely, and for
$P_\alpha=\alpha P_1+(1-\alpha)P_2$,
\[
 d_{P_\alpha}^0
 =
 \alpha d_{P_1}^0+(1-\alpha)d_{P_2}^0.
\]
Thus the construction embeds an infinite-dimensional convex family of
probability measures into the class of unit deviances generating PDMs.
This gives a stronger structural interpretation of the abundance of
PDMs than the cardinality statement alone.

\noindent

\subsection{Construction of Non-Standard Dispersion Models} \label{Sect.3.3}

We now show that a proper dispersion model of location type may, under
suitable conditions, be modified by replacing its constant normalising
factor by a non-constant positive function, yielding a non-standard dispersion model (NSDM). 
The construction is based on zeros of the Fourier transform of the exponential kernel associated
with the unit deviance.

Let
$d(y;\mu)=d_0(y-\mu)$
be a unit deviance such that, for every $\lambda>0$,
\[
        K_\lambda(t)
        :=
        \exp\{-\lambda d_0(t)\}
        \in L^1(\mathbb{R}).
\]
Since $d_0$ is even, $K_\lambda$ is even. The corresponding proper
dispersion model has density
\[
        p_0(y;\mu,\lambda)
        =
        \widetilde a_\lambda
        \exp\{-\lambda d_0(y-\mu)\},
\]
where
\begin{equation}
\label{eq:atilde}
        \widetilde a_\lambda
        =
        \left[
        \int_{\mathbb{R}}K_\lambda(t)\,dt
        \right]^{-1}.
\end{equation}

We use throughout this section the Fourier transform convention
\[
        \widehat K_\lambda(\omega)
        =
        \int_{\mathbb{R}}
        e^{-i\omega t}K_\lambda(t)\,dt.
\]
Since $K_\lambda$ is real and even,
$\widehat K_\lambda$ is also real and even.

The following proposition gives a simple mechanism for constructing a
non-constant normalising function.
\begin{proposition}
\label{prop:FourierZeroNSDM}
Suppose that, for some $\lambda_0>0$, there exists
$\omega_0>0$ such that
\begin{equation}
\label{eq:FourierZero}
        \widehat K_{\lambda_0}(\omega_0)=0.
\end{equation}
Then, for every $\varepsilon$ satisfying
$0<|\varepsilon|<\widetilde a_{\lambda_0}$,
the function
\begin{equation}
\label{eq:nonconstantnormalizer}
        a_\varepsilon(y;\lambda_0)
        =
        \widetilde a_{\lambda_0}
        +
        \varepsilon\cos(\omega_0y)
\end{equation}
is positive and satisfies
\begin{equation}
\label{eq:NSDMnormalization}
        \int_{\mathbb{R}}
        a_\varepsilon(y;\lambda_0)
        \exp\{-\lambda_0d_0(y-\mu)\}\,dy
        =
        1
\end{equation}
for every $\mu\in\mathbb{R}$.
Consequently, defining
\[
 a_\varepsilon(y;\lambda)
 =
 \begin{cases}
 \widetilde a_{\lambda_0}
 +\varepsilon\cos(\omega_0y),
       & \lambda=\lambda_0,\\[2mm]
 \widetilde a_\lambda,
       & \lambda\neq\lambda_0,
 \end{cases}
\]
gives a dispersion model on
$\mathbb{R}\times\mathbb{R}_+$ having a non-constant normalising
function.
\end{proposition}

\begin{proof}
Since
$|\cos(\omega_0y)|\leq1$,
we have
$a_\varepsilon(y;\lambda_0)
        \geq
        \widetilde a_{\lambda_0}-|\varepsilon|
        >0$.
We now verify the normalisation condition.

By the change of variable $t=y-\mu$,
\[
 \int_{\mathbb{R}}
 \cos(\omega_0y)
 K_{\lambda_0}(y-\mu)\,dy
 =
 \int_{\mathbb{R}}
 \cos\{\omega_0(t+\mu)\}
 K_{\lambda_0}(t)\,dt.
\]
Using
$\cos\{\omega_0(t+\mu)\}
 =
 \cos(\omega_0t)\cos(\omega_0\mu)
 -
 \sin(\omega_0t)\sin(\omega_0\mu)
 $,
and recalling that $K_{\lambda_0}$ is even, we obtain
$
 \int_{\mathbb{R}}
 \sin(\omega_0t)K_{\lambda_0}(t)\,dt
 =0
 $,
whereas
$
 \int_{\mathbb{R}}
 \cos(\omega_0t)K_{\lambda_0}(t)\,dt
 =
 \widehat K_{\lambda_0}(\omega_0)
$.
Thus, by \eqref{eq:FourierZero},
$
        \int_{\mathbb{R}}
        \cos(\omega_0y)
        K_{\lambda_0}(y-\mu)\,dy
        =0
$.
It follows that
\[
 \int_{\mathbb{R}}
 a_\varepsilon(y;\lambda_0)
 K_{\lambda_0}(y-\mu)\,dy
  =
 \widetilde a_{\lambda_0}
 \int_{\mathbb{R}}K_{\lambda_0}(y-\mu)\,dy
 =1,
\]
which proves \eqref{eq:NSDMnormalization}.
\end{proof}

\noindent
More generally, if
$\widehat K_{\lambda_0}(\omega_0)=0$, then both
$\cos(\omega_0y)$ and $\sin(\omega_0y)$ have zero integral against
every translate of $K_{\lambda_0}$. Thus one may use perturbations of
the form
$
        A\cos(\omega_0y)+B\sin(\omega_0y)
$,
provided their amplitude is sufficiently small to preserve positivity
of the normalising function.

The condition in Proposition \ref{prop:FourierZeroNSDM} is not
vacuous. We now give an explicit example of a unit deviance for which
the corresponding kernel has a real Fourier zero.

\begin{example}
\label{ex:FourierZeroKernel}
Let $\eta>0$ and $\beta>0$ satisfy
\begin{equation}
\label{eq:etabeta}
        \eta\beta^2<2,
\end{equation}
and define
\begin{equation}
\label{eq:explicitK}
        K(t)
        =
        e^{-t^2}
        \left[
        1+\eta\{1-\cos(\beta t)\}
        \right].
\end{equation}
Then
$K(0)=1$
and
$0<K(t)<1$, for $t\neq0$.
Indeed, positivity is immediate. Moreover,
\[
        1-\cos(\beta t)
        \leq \frac{\beta^2t^2}{2},
\]
and hence
\[
        1+\eta\{1-\cos(\beta t)\}
        \leq
        1+\frac{\eta\beta^2}{2}t^2.
\]
By \eqref{eq:etabeta},
\[
        \frac{\eta\beta^2}{2}<1,
\]
and therefore, for $t\neq0$,
\[
        1+\frac{\eta\beta^2}{2}t^2
        <
        e^{t^2}.
\]
Consequently,
\[
        K(t)<1
        \qquad (t\neq0).
\]

Define
\begin{equation}
\label{eq:explicitd0}
        d_0(t)
        =
        -\log K(t)
        =
        t^2
        -
        \log\left[
        1+\eta\{1-\cos(\beta t)\}
        \right].
\end{equation}
Then
\[
        d_0(0)=0,
        \qquad
        d_0(t)>0
        \quad(t\neq0),
\]
and therefore
\[
        d(y;\mu)=d_0(y-\mu)
\]
is a unit deviance.

Furthermore, since
\[
        1
        \leq
        1+\eta\{1-\cos(\beta t)\}
        \leq
        1+2\eta,
\]
we have
\[
        t^2-\log(1+2\eta)
        \leq d_0(t)\leq t^2.
\]
Thus
\[
        d_0(t)\longrightarrow\infty
        \qquad\text{as }|t|\longrightarrow\infty,
\]
and, for every $\lambda>0$,
\[
        e^{-\lambda d_0(t)}
\]
is integrable on $\mathbb{R}$.

For $\lambda_0=1$ the kernel associated with the deviance is precisely
the function $K$ in \eqref{eq:explicitK}. Its Fourier transform is
\begin{equation}
\label{eq:explicitKhat}
\widehat K(\omega)
=
\sqrt{\pi}
\left[
(1+\eta)e^{-\omega^2/4}
-
\frac{\eta}{2}
e^{-(\omega-\beta)^2/4}
-
\frac{\eta}{2}
e^{-(\omega+\beta)^2/4}
\right].
\end{equation}
At the origin,
\[
\widehat K(0)
=
\sqrt{\pi}
\left[
1+\eta-\eta e^{-\beta^2/4}
\right]
>0.
\]
On the other hand,
\[
\frac{
e^{-(\omega-\beta)^2/4}
}{
e^{-\omega^2/4}
}
=
\exp\left\{
\frac{\beta\omega}{2}
-
\frac{\beta^2}{4}
\right\}
\longrightarrow\infty
\qquad
(\omega\longrightarrow\infty).
\]
Hence the negative term
\[
        -\frac{\eta}{2}
        e^{-(\omega-\beta)^2/4}
\]
eventually dominates the positive term
\[
        (1+\eta)e^{-\omega^2/4},
\]
and therefore
\[
        \widehat K(\omega)<0
\]
for all sufficiently large positive $\omega$.

Since $\widehat K$ is continuous and
$\widehat K(0)>0$, there exists at least one
$\omega_0>0$ such that
\[
        \widehat K(\omega_0)=0.
\]
Proposition \ref{prop:FourierZeroNSDM} can therefore be applied at
$\lambda_0=1$.
\end{example}

For the unit deviance in Example \ref{ex:FourierZeroKernel}, choose any
zero $\omega_0>0$ of $\widehat K$ and any
\[
        0<|\varepsilon|<\widetilde a_1.
\]
Then
\[
        a_\varepsilon(y;1)
        =
        \widetilde a_1
        +
        \varepsilon\cos(\omega_0y)
\]
is a positive non-constant normalising function.

The corresponding dispersion model is not proper. Indeed, suppose that
its normalising function could be written in the form
\[
        a_\varepsilon(y;\lambda)
        =
        a_0(\lambda)b(y).
\]
For any $\lambda\neq1$ the function
$a_\varepsilon(y;\lambda)=\widetilde a_\lambda$ is constant in $y$.
It follows that $b(y)$ must be constant. But then
$a_\varepsilon(y;1)$ would also be constant in $y$, contradicting
$\varepsilon\neq0$.

Moreover, the unit deviance
\[
        d_0(t)
        =
        t^2
        -
        \log\left[
        1+\eta\{1-\cos(\beta t)\}
        \right]
\]
is non-quadratic whenever $\eta>0$ and $\beta>0$. Hence, by the
characterisation of location-type exponential dispersion models given
in the preceding section, the associated dispersion model is not an
EDM. It is therefore a non-standard dispersion model.

We record this conclusion explicitly.

\begin{corollary}
\label{cor:ExplicitNSDM}
There exist non-standard dispersion models on
$\mathbb{R}\times\mathbb{R}_+$ having positive non-constant
normalising functions. In particular, each choice of
$\eta,\beta>0$ satisfying $\eta\beta^2<2$, together with a real zero
$\omega_0$ of the Fourier transform in
\eqref{eq:explicitKhat}, gives the family
\[
        a_\varepsilon(y;1)
        =
        \widetilde a_1
        +
        \varepsilon\cos(\omega_0y),
        \qquad
        0<|\varepsilon|<\widetilde a_1.
\]
Thus a single unit deviance of the form
\eqref{eq:explicitd0} generates a family of distinct non-standard
dispersion models having cardinality at least that of the real line.
\end{corollary}

The construction above does not require the perturbing function to
belong to $L^2(\mathbb{R})$. The only essential requirements are that
the perturbation be bounded, so that positivity of the normalising
function can be preserved, and that its integral against every
translate of $K_{\lambda_0}$ vanishes. A zero of the Fourier transform
provides such a perturbation directly.

An alternative construction based on
$L^2(\mathbb{R})$-orthogonality is considered below. In that case the
zero set of $\widehat K_\lambda$ must have positive Lebesgue measure,
rather than merely containing a single point.

\subsection{An $L^2$ Construction of Non-Standard Dispersion Models} \label{Sect.3.4}

We consider now an alternative construction of non-standard dispersion
models based on orthogonality in $L^2(\mathbb{R})$. 
Throughout this subsection, $\nu$ denotes the Lebesgue measure on
$\mathbb{R}$.

Let $d(y;\mu)=d_0(y-\mu)$
be a unit deviance such that, for each $\lambda>0$,
\[
        K_\lambda(t)
        =
        \exp\{-\lambda d_0(t)\}
        \in L^1(\mathbb{R})\cap L^2(\mathbb{R}).
\]
As before, the corresponding proper dispersion model has constant
normalising function
\[
        \widetilde a_\lambda
        =
        \left[
        \int_{\mathbb{R}}K_\lambda(t)\,dt
        \right]^{-1}.
\]
For a fixed $\lambda>0$, consider the closed linear subspace
\begin{equation}
\label{eq:Vlambda}
 V_\lambda
 =
 \overline{\operatorname{span}
 \left\{
 K_\lambda(\cdot-\mu):\mu\in\mathbb{R}
 \right\}}^{\,L^2(\mathbb{R})}.
\end{equation}
A function $f\in L^2(\mathbb{R})$ belongs to
$V_\lambda^\perp$ if ,and only if,
\begin{equation}
\label{eq:L2orthogonality}
 \int_{\mathbb{R}}
 f(y)K_\lambda(y-\mu)\,dy
 =
 0,
 \qquad \forall \mu\in\mathbb{R}.
\end{equation}
Thus every suitably bounded non-zero element of
$V_\lambda^\perp$ can potentially be added, after multiplication by a
sufficiently small constant, to the trivial normalising function
$\widetilde a_\lambda$.

The structure of $V_\lambda^\perp$ can be characterised completely in
terms of the Fourier transform of $K_\lambda$. Define
\begin{equation}
\label{eq:Zlambda}
 Z_\lambda
 =
 \left\{
 \omega\in\mathbb{R}:
 \widehat K_\lambda(\omega)=0
 \right\}.
\end{equation}

\begin{proposition}
\label{prop:L2orthogonal}
For each $\lambda>0$,
\begin{equation}
\label{eq:VperpFourier}
 V_\lambda^\perp
 =
 \left\{
 f\in L^2(\mathbb{R}):
 \widehat f(\omega)=0
 \text{ for almost every }
 \omega\notin Z_\lambda
 \right\}.
\end{equation}
Consequently,
\begin{equation}
\label{eq:VperpCondition}
 V_\lambda^\perp\neq\{0\}
 \quad\Longleftrightarrow\quad
 \nu(Z_\lambda)>0.
\end{equation}
\end{proposition}

\begin{proof}
Let $f\in L^2(\mathbb{R})$. By Plancherel's theorem,
for every $\mu\in\mathbb{R}$,
\[
\begin{split}
 \left\langle
 f,K_\lambda(\cdot-\mu)
 \right\rangle_{L^2}
 &=
 \frac{1}{2\pi}
 \int_{\mathbb{R}}
 \widehat f(\omega)
 \overline{\widehat K_\lambda(\omega)}
 e^{i\omega\mu}\,d\omega .
\end{split}
\]
Here
$
 \widehat f\,
 \overline{\widehat K_\lambda}
 \in L^1(\mathbb{R})
$
by the Cauchy--Schwarz inequality. 
Therefore,
$
 f\in V_\lambda^\perp
$
if, and only if, the Fourier transform of the $L^1$ function
$
 \widehat f\,
 \overline{\widehat K_\lambda}
$
vanishes identically. By uniqueness of the Fourier transform,
this is equivalent to
$
 \widehat f(\omega)
 \overline{\widehat K_\lambda(\omega)}
 =
 0
 $for almost every
 $\omega\in\mathbb{R}$.
Hence,
$ \widehat f(\omega)=0$
for almost every $\omega\notin Z_\lambda $,
which proves \eqref{eq:VperpFourier}.

If $\nu(Z_\lambda)=0$, any $L^2$ function supported, up to sets of
measure zero, on $Z_\lambda$ is equal to zero almost everywhere, and
therefore $V_\lambda^\perp=\{0\}$.

Conversely, suppose that $\nu(Z_\lambda)>0$. Since the Lebesgue measure is
$\sigma$-finite, there exists a measurable set
$E\subseteq Z_\lambda$ such that
$0<\nu(E)<\infty$.
Then $\mathbf{1}_E\in L^2(\mathbb{R})$ is non-zero, and its inverse
Fourier transform belongs to $V_\lambda^\perp$. Hence
$V_\lambda^\perp\neq\{0\}$.
\end{proof}

\noindent
The preceding result also provides bounded elements of
$V_\lambda^\perp$, which are useful for the construction of positive
normalising functions.

\begin{proposition}
\label{prop:boundedL2perturbation}
Suppose that
$\nu(Z_{\lambda_0})>0$
for some $\lambda_0>0$. 
Then there exists a non-zero real-valued function
$
        f\in V_{\lambda_0}^\perp
        \cap L^\infty(\mathbb{R})
$.
Consequently, for every $\varepsilon$ satisfying
\begin{equation}
\label{eq:L2epsilon}
 0<|\varepsilon|
 <
 \frac{\widetilde a_{\lambda_0}}
      {\|f\|_\infty},
\end{equation}
the function
\begin{equation}
\label{eq:L2normalizer}
 a_\varepsilon(y;\lambda_0)
 =
 \widetilde a_{\lambda_0}
 +
 \varepsilon f(y)
\end{equation}
is positive and satisfies
\[
 \int_{\mathbb{R}}
 a_\varepsilon(y;\lambda_0)
 K_{\lambda_0}(y-\mu)\,dy
 =
 1
\]
for every $\mu\in\mathbb{R}$.
\end{proposition}

\begin{proof}
Since $K_{\lambda_0}$ is real and even,
$\widehat K_{\lambda_0}$ is real and even, and hence
$Z_{\lambda_0}$ is symmetric about zero.
Since $\nu(Z_{\lambda_0})>0$, we may choose a symmetric measurable
set
$
        E\subseteq Z_{\lambda_0}
$
such that
$
        0<\nu(E)<\infty
$.
Define $f$ by
\begin{equation}
\label{eq:fviaE}
 \widehat f(\omega)=\mathbf{1}_E(\omega).
\end{equation}
Since
$
        \mathbf{1}_E\in L^1(\mathbb{R})
        \cap L^2(\mathbb{R}),
$
Fourier inversion and Plancherel's theorem imply that
$
        f\in L^\infty(\mathbb{R})
        \cap L^2(\mathbb{R})
$.
Moreover, $f$ is continuous, real-valued and even. Since
$\mathbf{1}_E$ is supported in $Z_{\lambda_0}$,
Proposition \ref{prop:L2orthogonal} gives that
$
        f\in V_{\lambda_0}^\perp
$.

Condition \eqref{eq:L2epsilon} implies
$
 a_\varepsilon(y;\lambda_0)
 \geq
 \widetilde a_{\lambda_0}
 -
 |\varepsilon|\,\|f\|_\infty\\
 >0
$
for every $y\in\mathbb{R}$. Finally,
because $f\in V_{\lambda_0}^\perp$,
\[
 \int_{\mathbb{R}}
 a_\varepsilon(y;\lambda_0)
 K_{\lambda_0}(y-\mu)\,dy
 =
 \widetilde a_{\lambda_0}
 \int_{\mathbb{R}}K_{\lambda_0}(y-\mu)\,dy
 +
 \varepsilon
 \int_{\mathbb{R}}
 f(y)K_{\lambda_0}(y-\mu)\,dy
 =1.
\]
\end{proof}

As in the construction of Subsection \ref{Sect.3.3}, one may define
\[
 a_\varepsilon(y;\lambda)
 =
 \begin{cases}
 \widetilde a_{\lambda_0}+\varepsilon f(y),
       &\lambda=\lambda_0,\\[2mm]
 \widetilde a_\lambda,
       &\lambda\neq\lambda_0.
 \end{cases}
\]
This gives a dispersion model with a non-constant normalising
function. The model cannot be proper: if
$
        a_\varepsilon(y;\lambda)=a_0(\lambda)b(y)
$,
then the fact that $a_\varepsilon(y;\lambda)$ is constant in $y$ for
every $\lambda\neq\lambda_0$ implies that $b$ is constant, contrary to
the non-constancy of $a_\varepsilon(\cdot;\lambda_0)$.
If, in addition, the unit deviance $d$ is not the unit deviance of an
exponential dispersion model, the resulting model is therefore a
non-standard dispersion model.

Condition \eqref{eq:VperpCondition} shows precisely the difference
between the present construction and that of Subsection
\ref{Sect.3.3}. A single real zero of $\widehat K_\lambda$ is
sufficient to construct a bounded trigonometric perturbation as in
Subsection \ref{Sect.3.3}. A non-zero perturbation belonging to
$L^2(\mathbb{R})$, however, exists if and only if the zero set of
$\widehat K_\lambda$ has positive Lebesgue measure.

\subsubsection*{An explicit example}

We conclude by constructing explicitly a unit deviance for which the
orthogonal complement $V_\lambda^\perp$ is infinite-dimensional for
at least one value of the index parameter.

Let $\{U_n:n\geq1\}$ be a sequence of independent random variables
such that
$
        U_n\sim {\rm Unif}[-2^{-n},2^{-n}]
$,
and define
$
        X=\sum_{n=1}^\infty U_n 
$.
Since
$
        \sum_{n=1}^\infty 2^{-n}=1
$,
the random variable $X$ is well defined and has support contained in
$[-1,1]$. Its characteristic function is
\begin{equation}
\label{eq:InfiniteProductCF}
 \varphi(t)
 =
 \prod_{n=1}^\infty
 \frac{\sin(t/2^n)}{t/2^n}.
\end{equation}

The distribution of $X$ has a symmetric $C^\infty$ density $q$ with
compact support. Indeed, the product in
\eqref{eq:InfiniteProductCF} decreases faster than every negative
power of $|t|$. To see this, use
$
 \left|
 \frac{\sin x}{x}
 \right|
 \leq
 \min\{1,|x|^{-1}\}
$
on the initial factors of the product. It follows that, for every
$m>0$,
$
        |t|^m|\varphi(t)|\longrightarrow0
         (|t|\longrightarrow\infty)
$.
Fourier inversion then implies that the distribution of $X$ possesses
a $C^\infty$ density $q$, while the construction of $X$ implies that
$q$ is supported on $[-1,1]$.

The non-zero real zeros of $\varphi$ are particularly simple. The
$n$th factor in \eqref{eq:InfiniteProductCF} vanishes when
$ t=k\pi 2^n$, $k\in\mathbb Z\setminus\{0\}$.
The union of these sets, over $n\geq1$, is
$2\pi\mathbb Z\setminus\{0\}$.
There are no further zeros. Indeed, if no factor vanishes, then
\[
 \sum_{n=1}^\infty
 \left|
 1-
 \frac{\sin(t/2^n)}{t/2^n}
 \right|
 <\infty,
\]
and hence the infinite product is non-zero. Thus
\begin{equation}
\label{eq:zerosphi}
 \{t\in\mathbb R:\varphi(t)=0\}
 =
 2\pi\mathbb Z\setminus\{0\}.
\end{equation}

Take now
$
        r=\sqrt{2}
$.
If $\varphi(t)=\varphi(rt)=0$ for some $t\neq0$, then by
\eqref{eq:zerosphi} there would exist non-zero integers $k$ and $l$
such that
$t=2\pi k$ and $\sqrt{2}\,t=2\pi l$.
This would imply
$
        \sqrt{2}=\frac{l}{k}
$,
which is impossible. Hence $\varphi(t)$ and
$\varphi(\sqrt{2}\,t)$ have no common non-zero real zero.

Define
\begin{equation}
\label{eq:compactSpectrumK}
 K(t)
 =
 \frac{1}{2}
 \left\{
 \varphi(t)^2+\varphi(\sqrt{2}\,t)^2
 \right\}.
\end{equation}
Then
$
        K(0)=1
$,
and
$
        K(t)>0, t\neq0
$,
because the two terms in \eqref{eq:compactSpectrumK} cannot vanish
simultaneously.

Furthermore,
$
        K(t)<1, t\neq0
$.
Indeed, since the distribution of $X$ is absolutely continuous,
$
        |\varphi(t)|<1, t\neq0
$,
and the same holds for $\varphi(\sqrt{2}\,t)$.

It follows that
\begin{equation}
\label{eq:compactSpectrumDeviance}
        d_0(t)=-\log K(t)
\end{equation}
satisfies
$d_0(0)=0$ and $d_0(t)>0, t\neq0$.
Therefore
$d(y;\mu)=d_0(y-\mu)$
is a unit deviance.

Since $\varphi$ decreases faster than every negative power,
the same is true of $K$. Consequently, for every $\lambda>0$ and
every sufficiently large $N$,
\[
        K(t)^\lambda
        \leq
        C_N(1+|t|)^{-N},
\]
and hence
$
        K_\lambda(t)
        =
        \exp\{-\lambda d_0(t)\}
        =
        K(t)^\lambda
$
belongs to $L^1(\mathbb R)\cap L^2(\mathbb R)$. Thus $d$ generates a
proper dispersion model for every $\lambda>0$.

We now consider $\lambda=1$. Let $q*q$ denote the convolution of
$q$ with itself. With the Fourier transform convention used in this
paper,
\[
        \widehat{\varphi^2}(\omega)
        =
        2\pi(q*q)(\omega),
\]
whereas
\[
        \widehat{\varphi(\sqrt{2}\,\cdot)^2}(\omega)
        =
        \frac{2\pi}{\sqrt{2}}
        (q*q)\left(\frac{\omega}{\sqrt{2}}\right).
\]
Hence
\begin{equation}
\label{eq:KhatCompactExample}
 \widehat K(\omega)
 =
 \pi
 \left\{
 (q*q)(\omega)
 +
 \frac{1}{\sqrt{2}}
 (q*q)\left(\frac{\omega}{\sqrt{2}}\right)
 \right\}.
\end{equation}

Since $q$ is supported on $[-1,1]$, the convolution $q*q$ is
supported on $[-2,2]$. It follows from
\eqref{eq:KhatCompactExample} that
$ \widehat K(\omega)=0$ for $ |\omega|>2\sqrt{2}$.
Therefore, with $Z_1=\{\omega:\widehat K(\omega)=0\}$,
we have
 $(-\infty,-2\sqrt{2})\cup(2\sqrt{2},\infty)
        \subseteq Z_1
$,
and consequently
$
        \nu(Z_1)=\infty
$.
By Proposition \ref{prop:L2orthogonal}, the space
$V_1^\perp$ is infinite-dimensional.

For example, choose
$
 E=[A,B]\cup[-B,-A]
 2\sqrt{2}<A<B,
$
and define
$
        \widehat f=\mathbf{1}_E
$.
Then
\[
 f(y)
 =
 \frac{\sin(By)-\sin(Ay)}{\pi y},
\]
with the continuous extension
$
        f(0)=\frac{B-A}{\pi}
$.
Thus
$
        f\in
        V_1^\perp
        \cap L^2(\mathbb R)
        \cap L^\infty(\mathbb R)
$,
and
$
        \widetilde a_1+\varepsilon f
$
is a positive non-constant normalising function whenever
$
        0<|\varepsilon|
        <
        \{\widetilde a_1\}/\|f\|_\infty
$.

Finally, the unit deviance
\eqref{eq:compactSpectrumDeviance} is not quadratic. If it were of the
form
$
        d_0(t)=ct^2
$,
then
$
        K(t)=e^{-ct^2}
$,
whose Fourier transform is a Gaussian and therefore has full support
on $\mathbb R$. This contradicts
\eqref{eq:KhatCompactExample}. Hence the unit deviance is not the unit
deviance of an exponential dispersion model, and the models obtained
from the non-constant normalising functions above are non-standard
dispersion models.

\section{Discussion and Examples}
\label{Sect.4}

The constructions developed in the preceding sections illustrate that
the class of dispersion models is considerably larger than the
classical families usually considered in applications. In particular,
characteristic functions provide a large source of unit deviances and
proper dispersion models, while zeros of the Fourier transform of the
associated exponential kernels may be used to construct non-constant
normalising functions and hence non-standard dispersion models.

A unit deviance $d$ is said to be \emph{regular} when it is twice
continuously differentiable and
\[
        \frac{\partial^2 d}{\partial\mu^2}(\mu;\mu)>0,
        \qquad \mu\in\Omega .
\]
If $\yyy=\Omega$, the corresponding dispersion model is called a
\emph{regular dispersion model}. Regularity plays an important role in
the theory of statistical inference for dispersion models
\citep{Jorgensen1997,Cordeiroetal2019}, since the second-order local
behaviour of the deviance determines, among other quantities, the unit
variance function.

The construction of Subsection \ref{Sect.3.2} preserves the
local behaviour of the characteristic-function component of the unit
deviance. Indeed, if
\[
        d_{\varphi,a,b,c}(y;\mu)
        =
        1-\varphi(y-\mu)
        +
        c(y-\mu)^2\rho_{a,b}(y-\mu),
\]
then, for $|y-\mu|\leq a$,
\[
        d_{\varphi,a,b,c}(y;\mu)
        =
        1-\varphi(y-\mu).
\]
Consequently, if the probability distribution associated with
$\varphi$ has finite and non-zero variance $\sigma^2$, then
\[
        1-\varphi(t)
        =
        \frac{\sigma^2}{2}t^2+o(t^2),
        \qquad t\longrightarrow0,
\]
and hence
\[
        \frac{\partial^2 d_{\varphi,a,b,c}}
             {\partial\mu^2}(\mu;\mu)
        =
        \sigma^2>0.
\]
Provided that the characteristic function and the cut-off function
have the required smoothness, the resulting dispersion model is
therefore regular.

\begin{figure}
  \centering
  \includegraphics[width=0.95\textwidth]{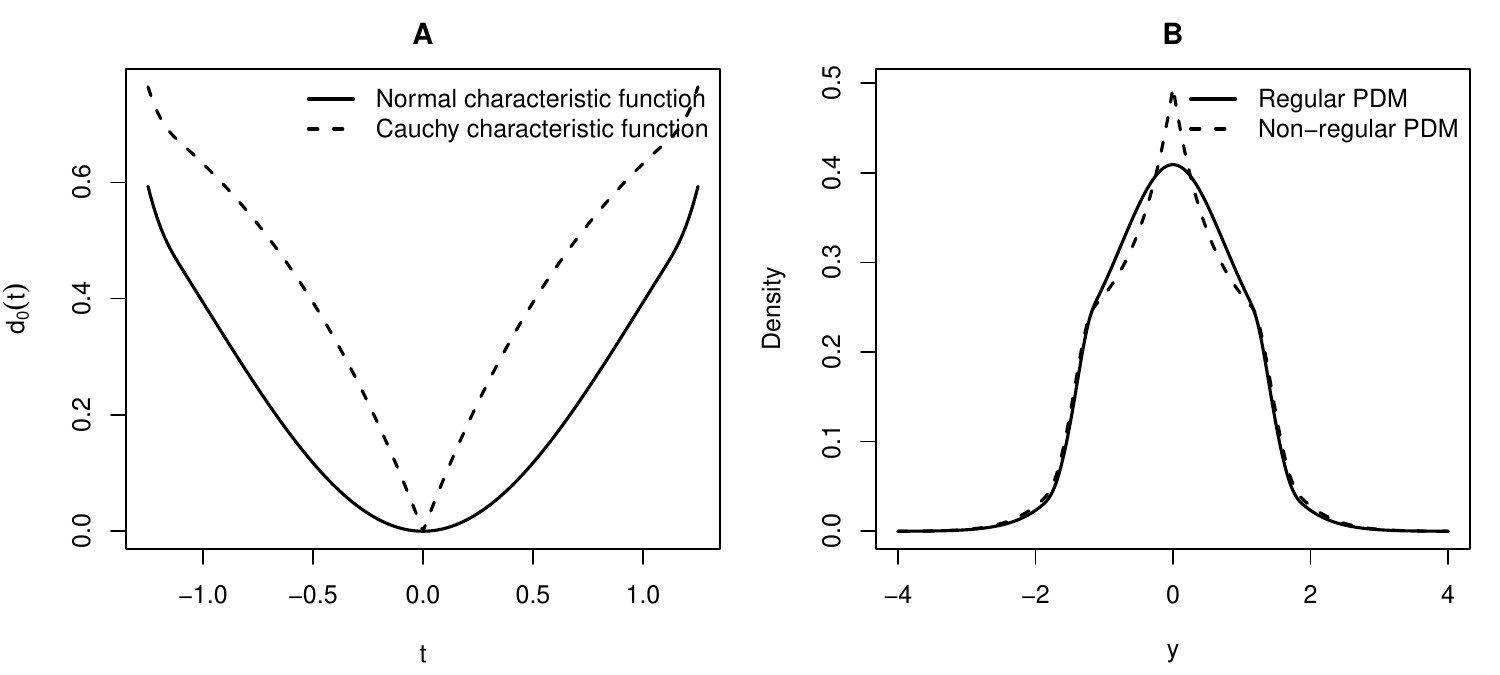}
  \caption{
  Examples of proper dispersion models generated by the construction
  of Subsection \ref{Sect.3.2}. Panel A shows the unit
  deviances obtained from the characteristic functions of the standard
  normal distribution (solid line) and the standard Cauchy distribution
  (dashed line), using $a=1$, $b=2$ and $c=1/2$. The normal-based
  unit deviance is regular, whereas the Cauchy-based unit deviance has
  a cusp at the origin and is non-regular. Panel B shows the
  corresponding densities for $\mu=0$ and $\lambda=1$.
  }
  \label{figure-PDM}
\end{figure}

On the other hand, the same construction also produces non-regular
dispersion models. For example, consider respectively the
characteristic functions of the standard normal and standard Cauchy
distributions,
\[
        \varphi_N(t)=\exp(-t^2/2)
        \qquad\hbox{and}\qquad
        \varphi_C(t)=\exp(-|t|).
\]
Taking, for example,
\[
        a=1,\qquad b=2,\qquad c=\frac12,
\]
and a smooth cut-off function $\rho_{a,b}$ as in Subsection
\ref{Sect.3.2}, we obtain the unit deviances
\[
 d_N^0(t)
 =
 1-\exp(-t^2/2)
 +
 \frac12t^2\rho_{1,2}(t)
\]
and
\[
 d_C^0(t)
 =
 1-\exp(-|t|)
 +
 \frac12t^2\rho_{1,2}(t).
\]
The first is regular, whereas the second is not differentiable at the
origin and hence is not regular. Nevertheless, both generate proper
dispersion models for every $\lambda>0$. Figure
\ref{figure-PDM} illustrates the local behaviour of the two unit
deviances and the corresponding densities for $\lambda=1$ and
$\mu=0$.

A different phenomenon is illustrated by the construction of
Subsection \ref{Sect.3.3}. Consider the kernel
\[
 K(t)
 =
 e^{-t^2}
 \left[
 1+\eta\{1-\cos(\beta t)\}
 \right],
\]
where $\eta,\beta>0$ and $\eta\beta^2<2$. The corresponding unit
deviance is
\[
 d_0(t)
 =
 t^2
 -
 \log\left[
 1+\eta\{1-\cos(\beta t)\}
 \right].
\]
For $\lambda=1$ the Fourier transform of the kernel is
\[
\widehat K(\omega)
=
\sqrt{\pi}e^{-\omega^2/4}
\left[
(1+\eta)
-
\eta e^{-\beta^2/4}
\cosh\left(\frac{\beta\omega}{2}\right)
\right].
\]
Consequently, $\widehat K$ has a positive real zero given explicitly by
\[
 \omega_0
 =
 \frac{2}{\beta}
 \operatorname{arcosh}
 \left[
 \frac{1+\eta}{\eta}
 e^{\beta^2/4}
 \right].
\]
For example, taking $\eta=\beta=1$ gives
\[
        \omega_0\simeq 3.192.
\]

At this frequency,
\[
        \int_{\mathbb{R}}
        \cos(\omega_0y)K(y-\mu)\,dy=0
\]
for every $\mu$. Hence
\[
        a_\varepsilon(y;1)
        =
        \widetilde a_1
        +
        \varepsilon\cos(\omega_0y),
        \qquad
        0<|\varepsilon|<\widetilde a_1,
\]
is a positive non-constant normalising function. The resulting
dispersion model is not proper. Moreover, since $d_0$ is non-quadratic,
it is not an exponential dispersion model and hence is an NSDM.

\begin{figure}
  \centering
  \includegraphics[width=0.98\textwidth]{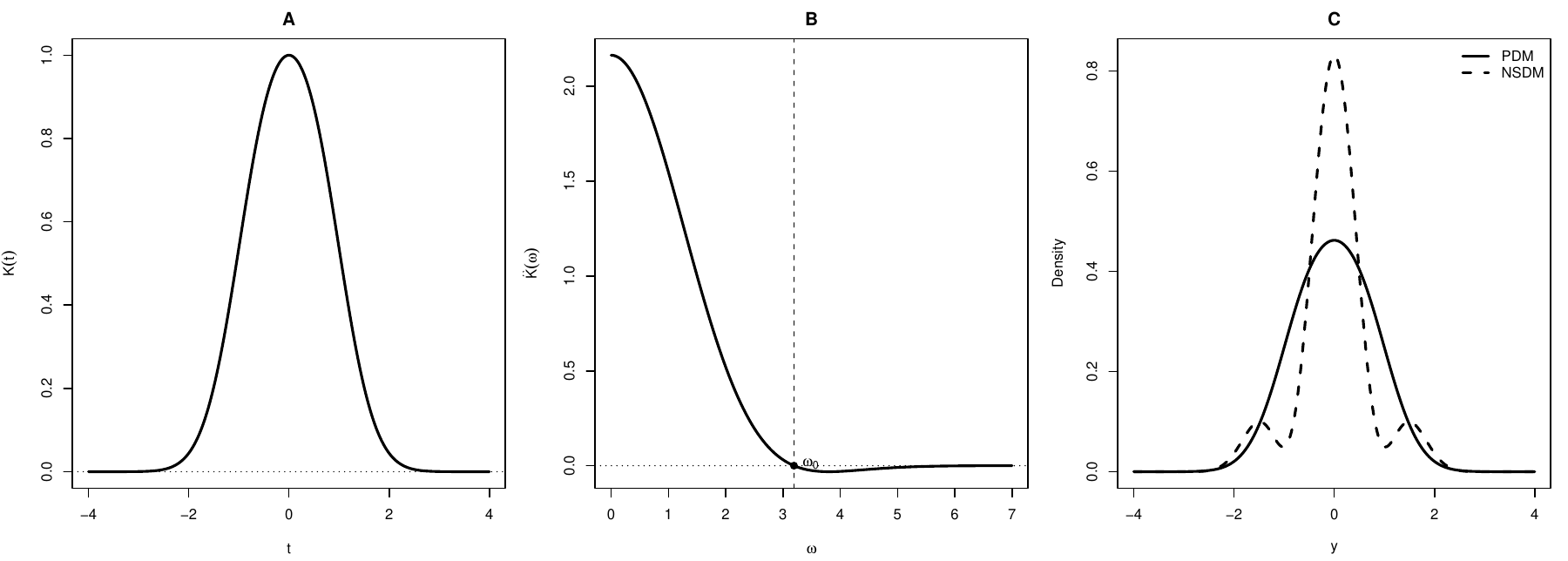}
  \caption{
  Construction of an NSDM using a single zero of the Fourier transform
  of the exponential kernel. Here $\eta=\beta=1$ and $\lambda=1$.
  Panel A shows
  $K(t)=e^{-t^2}\{1+\eta[1-\cos(\beta t)]\}$.
  Panel B shows $\widehat K$ and its first positive zero
  $\omega_0$. Panel C compares the density of the associated PDM
  (solid line) with an NSDM density obtained using
  $a_\varepsilon(y;1)
   =\widetilde a_1+\varepsilon\cos(\omega_0y)$
  (dashed line).
  }
  \label{figure-NSDM-Fourier}
\end{figure}

Figure \ref{figure-NSDM-Fourier} illustrates this construction for
$\eta=\beta=1$. Panel A shows the kernel $K$, panel B shows its Fourier
transform and the positive zero $\omega_0$, and panel C compares the
density of the corresponding PDM with a density obtained by replacing
the constant normalising function by
$a_\varepsilon(y;1)$. The oscillatory modification of the normalising
function may produce shapes which cannot arise in a PDM generated by
the same unit deviance.

The $L^2$ construction of Subsection \ref{Sect.3.3} provides a
different perspective. For a fixed $\lambda$, let
$
 Z_\lambda
 =
 \{\omega\in\mathbb{R}:
   \widehat K_\lambda(\omega)=0\}
$.
The space of $L^2(\mathbb{R})$ functions orthogonal to all translates
of $K_\lambda$ is non-trivial if and only if
$
        \nu(Z_\lambda)>0
$,
where, in this context, $\nu$ denotes Lebesgue measure on
$\mathbb{R}$. Thus the $L^2$ construction requires more than an
isolated Fourier zero: the zero set must have positive Lebesgue
measure.

The example given in Subsection \ref{Sect.3.4} makes this condition
particularly transparent. Let
\[
        X=\sum_{n=1}^\infty U_n,
        \qquad
        U_n\sim{\rm Unif}[-2^{-n},2^{-n}]
\]
independently. Its characteristic function is
\[
        \varphi(t)
        =
        \prod_{n=1}^\infty
        \frac{\sin(t/2^n)}{t/2^n},
\]
and its non-zero real zeros are precisely
$2\pi\mathbb Z\setminus\{0\}$. Consequently,
$\varphi(t)$ and $\varphi(\sqrt{2}\,t)$ have no common non-zero zero.
Defining
\[
 K(t)
 =
 \frac12
 \left\{
 \varphi(t)^2+\varphi(\sqrt{2}\,t)^2
 \right\},
\]
we obtain a strictly positive kernel whose Fourier transform is
supported on
$
        [-2\sqrt{2},2\sqrt{2}]
$.
Thus
$
        \nu(Z_1)=\infty
$
and $V_1^\perp$ is infinite-dimensional.

\begin{figure}
  \centering
  \includegraphics[width=0.98\textwidth]{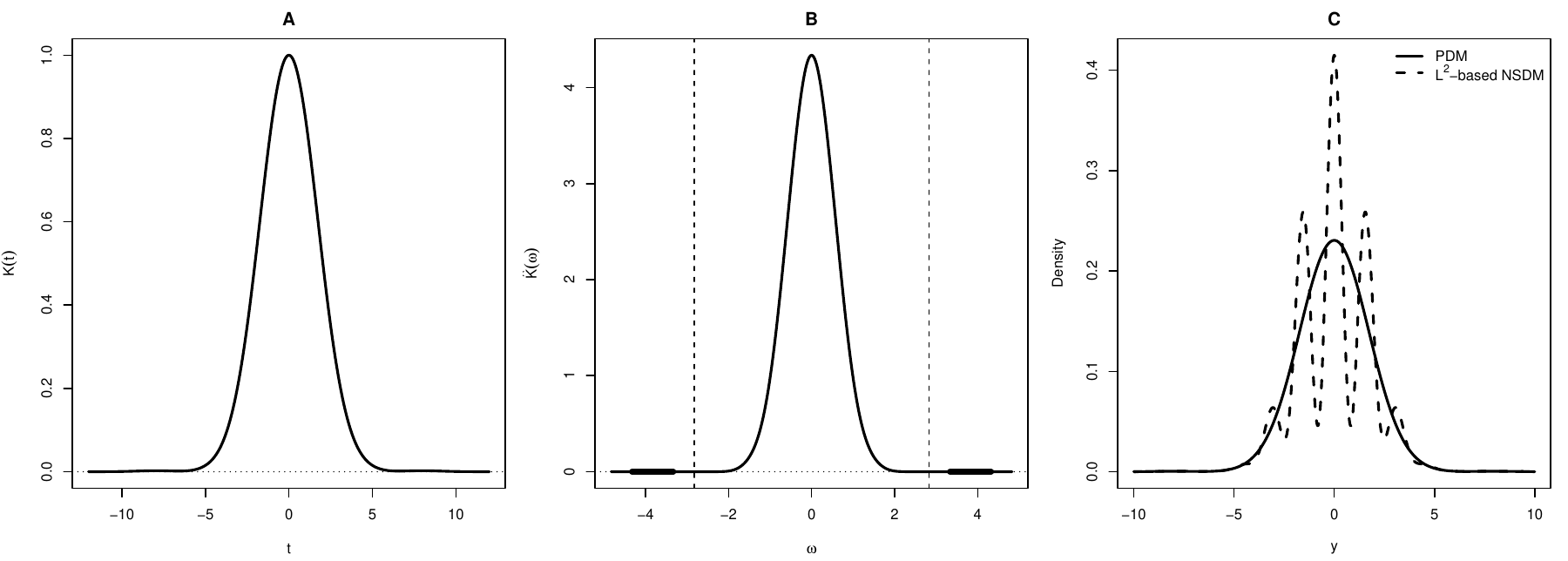}
  \caption{
  The $L^2$ construction of a non-standard dispersion model.
  Panel A shows
  $K(t)=\{\varphi(t)^2+\varphi(rt)^2\}/2$ for a characteristic
  function $\varphi$ obtained from a symmetric smooth density with
  compact support. Panel B shows the compactly supported Fourier
  transform $\widehat K$; the thick segments indicate a symmetric
  finite-measure set $E$ contained in the zero set of $\widehat K$.
  Panel C compares the associated PDM density (solid line) with an
  NSDM density obtained by adding a sufficiently small multiple of
  the inverse Fourier transform of $\mathbf{1}_E$ to the constant
  normalising function (dashed line).
  }
  \label{figure-NSDM-L2}
\end{figure}

Figure \ref{figure-NSDM-L2} illustrates this construction. Panel A
shows the kernel $K$, panel B its compactly supported Fourier transform
together with a set $E$ contained in the spectral zero set, and panel C
compares the original PDM density with an NSDM obtained by the
$L^2$-orthogonal perturbation described above.

The three examples above emphasise different aspects of the
construction. In Subsection \ref{Sect.3.2}, the characteristic
function determines the local geometry of the unit deviance, while the
additional coercive term controls the behaviour at infinity and
guarantees normalisability. In Subsection \ref{Sect.3.3}, a single real
zero of the Fourier transform of the exponential kernel is enough to
construct a bounded non-constant normalising function. In Subsection
\ref{Sect.3.4}, a zero set of positive Lebesgue measure yields an
entire infinite-dimensional space of $L^2$ perturbations.

These constructions also show that the distinction between proper,
exponential and non-standard dispersion models is structural rather
than merely terminological. A single unit deviance may generate a
proper dispersion model through its constant normalising function and,
when the associated kernel has suitable spectral zeros, additional
non-standard dispersion models through non-constant normalising
functions. The variety of such constructions indicates that the class
of dispersion models is substantially richer than the classical
parametric examples alone might suggest.

\section*{Acknowledgements}

We thank  Ole E. Barndorff-Nielsen (Department of Mathematics, Aarhus University), Denise A. Botter (Universidade de S\~ao Paulo) and Gauss M. Cordeiro (Universidade Federal de Pernambuco) for helpful comments in the early stage of this work.

\end{document}